\documentclass[11pt]{article}
\usepackage{amsmath}
\usepackage{amssymb}
\usepackage{theorem}

\usepackage{enumerate}

\usepackage{mathrsfs}

\usepackage{geometry}

\usepackage{hyperref}
\usepackage{graphicx}
\usepackage{caption}
\usepackage{subcaption}
\usepackage{epstopdf}

\numberwithin{equation}{section}

\newtheorem{thm}{Theorem}[section]
\newtheorem{lemma}[thm]{Lemma}
\newtheorem{prop}[thm]{Proposition}
\newtheorem{cor}[thm]{Corollary}
{\theorembodyfont{\rmfamily}

\newtheorem{rmk}[thm]{Remark}

}

\newcommand{\qed}{\hfill \mbox{\raggedright \rule{.07in}{.1in}}}

\newenvironment{proof}{\vspace{1ex}\noindent{\bf
Proof}\hspace{0.5em}}{\hfill\qed\vspace{1ex}}

\let\phi\varphi

\let\epsilon\varepsilon

\let\tilde\widetilde

\newcommand{\R}{{\mathbb R}}

\def\R{\ensuremath{\mathbb R}}

\def\dist{\ensuremath{\text{dist}}}

\begin{document}
\title{Dynamical Borel-Cantelli lemmas and rates of growth of Birkhoff sums of non-integrable observables on chaotic dynamical systems.}
\author{Meagan Carney \thanks{Department of Mathematics, University of Houston, Houston, USA. e-mail: $<$meagan@math.uh.edu$>$ 
Meagan Carney thanks the NSF for partial support on NSF-DMS Grant 1600780.}
\and
M. Nicol \thanks{Department of Mathematics,
University of Houston,
Houston Texas,
USA. e-mail: $<$nicol@math.uh.edu$>$.
 Matthew Nicol thanks the NSF for partial support on NSF-DMS Grant 1600780. Both authors thank Alan Haynes for helpful conversations.
}}

\date{\today}
\maketitle

\begin{abstract}
We consider implications of dynamical Borel-Cantelli lemmas for
rates of growth of Birkhoff sums of non-integrable observables $\phi(x)=d(x,p)^{-k}$, $k>0$, on ergodic dynamical systems $(T,X,\mu)$ where $\mu (X)=1$. Some general results are given as well as some more concrete examples involving  non-uniformly expanding maps, intermittent type maps as well as uniformly  hyperbolic systems.

\end{abstract}

\section{Birkhoff sums of non-integrable functions.}

Let $X_i$ be a sequence of random variables on a probability space $(X,\mu)$ (in other words a stochastic process) and let $S_n=\sum_{i=1}^n X_i$ be the associated sequence of 
Birkhoff sums.

W.~Feller~\cite{Feller} showed that if $\{X_i\}$ are iid and  $E|X_1|=\infty$ then for any sequence $b(n)>0$,  if $\lim_{n\to \infty} \frac{b(n)}{n}=\infty$ then either 
$\limsup \frac{S_n}{b(n)} =\infty $ a.e. or $\liminf \frac{S_n}{b(n)} =0 $ a.e. Chow and Robbins~\cite{Chow_Robbins} then showed that the conditions on $b(n)$ can be relaxed and that 
in fact for any sequence of constants $b(n)$ either $\limsup \frac{S_n}{b(n)} =\infty $ a.e. or $\liminf \frac{S_n}{b(n)} =0 $ a.e.

Suppose now that $(T, X, \mu)$ is an ergodic probability measure preserving transformation and $\phi : X\in \R$ is a non-integrable measurable function. Then $X_i:=\phi\circ T^i$
is a stationary stochastic process with  Birkhoff sum $S_n=\sum_{i=1}^n \phi \circ T^i$. In this dynamical setting 
Aaronson~\cite{Aaronson} showed that for any sequence $b(n)>0$, if $\lim_{n\to \infty} \frac{b(n)}{n}=\infty$ then either 
$\limsup \frac{S_n}{b(n)} =\infty $ a.e. or $\liminf \frac{S_n}{b(n)} =0 $ a.e. Thus for ergodic dynamical systems there is no strong law of large numbers for non-integrable observables.

A natural question is the rate of growth of Birkhoff sums. A useful result, due again to Aaronson~\cite[Proposition 2.3.1]{Aaronson} states:
\begin{prop}\label{upbound}
If $a(x)$ is increasing, $\lim_{x\to \infty} \frac{a(x)}{x}=0$ and
\[
\int a(\phi (x)) d\mu <\infty
\]
then for $\mu$ a.e. $x$ 
\[
\lim_{n\to \infty} \frac{a(S_n)}{n}=0
\]
\end{prop}
Despite the generality of its assumptions, if $p\in (0,1)$.
 gives   close to optimal bounds on $\limsup S_n$ in many dynamical settings, as demonstrated later in this paper.  Throughout this paper if $a(n)$ and $b(n)$ are two sequences $a(n)\sim b(n)$ will mean that there exists an $N$ and constants $C_1$, $C_2$  such that $0<C_1 \le \frac{a(n)}{b(n)} \le C_2$ for all $n\ge N$.

In~\cite{HNT} dynamical Borel-Cantelli lemmas were used to give information on the almost sure behavior of the maxima $M_n:=\max\{ \phi(x), \phi( T x),\phi (T^2 x), \ldots,  \phi(T^n x)\}$
for certain classes of observables on a variety of chaotic dynamical systems $(T,X,\mu)$. Motivated by applications in extreme value theory the observables considered in~\cite{HNT} were of  form 
$\phi(x)=-\log d(x,p)$ and  $\phi(x)=d(x,p)^{-k}$, where $d(.,.)$ was a Riemannian metric on the space $X$. For the integrable observable $\phi (x) =-\log d(x,p)$ under relatively mixing conditions on the dynamical system (please see~\cite[Theorem 2.2]{HNT} for details) a sequence of scaling constants $a(n)$ exists such that $\lim_{n\to \infty}\frac{M_n}{a(n)}=C>0$ almost surely for some constant $C$.

But the following result shows that for many dynamical systems  there is no almost sure limit for $\frac{M_n}{a(n)}$, in the case $\phi(x)=d(x,p)^{-k}$,  $k>0$, even if $k$  is such that $\phi$ is integrable (so that a  strong law of large numbers does hold for the Birkhoff sum). We state a simpler, less general,  version of~\cite[Theorem 2.7]{HNT} adapted for our purposes,

\begin{prop}
Suppose that $(T,X,\mu)$ is a measure preserving system with ergodic measure
$\mu$ which is absolutely continuous with respect to Lebesgue measure $m$. Suppose for 
a point $p\in X$ there exists $\delta>0$, $C>0$ and $r_0>0$ such that for all $\epsilon<r<r_0$:
\begin{equation}
|\mu(B(p,r+\epsilon))-\mu(B(p,r))|\leq C\epsilon^{\delta}.
\end{equation}
and  $0<\frac{d\mu}{dm}(p)<\infty$.
Moreover suppose that we have exponential decay of correlations in bounded variation norm (BV) versus $L^1$ in the sense that there exists $C>0$ and $0<\theta<1$ such that 
 for all $\varphi_1 $ of bounded variation and all $\varphi_2\in L^1 (m) $ we have:
\[
\left|\int \varphi_1 \cdot \varphi_2\circ f^j d \mu -\int \varphi_1 d \mu \int \varphi_2 d\mu\right|
\leq
C\theta^j \|\varphi_1\|_{BV} \|\varphi_2\|_{L^1(m)},
\]
Then if
 $\phi(x)=\dist(x,p)^{-k}$ for some $k>0$ for any monotone sequence $u(n)\to\infty$:
\begin{equation}\label{eq.mn.nolimit.thm}
\mu\left(\limsup_{n\to\infty}\frac{M_n(x)}{u(n)}=0\right)=1,\;\textrm{or}\;\mu\left(
\limsup_{n\to\infty}\frac{M_n(x)}{u(n)}=\infty\right)=1.
\end{equation}
\end{prop}

The relation between Birkhoff sums and extreme values, such as the maxima, is investigated in the topic of trimmed Birkhoff sums~\cite{Aaronson_Nakada,Kesten_Maller,Schindler}. In this approach
the time series $\{\phi (x), \phi (Tx), \phi (T^2 x),\ldots, \phi (T^n x)\}$ is rearranged into increasing order $\{\phi (T^{i_0} x) \le \phi (T^{i_1} x) \le  \phi (T^{i_2} x) \le \ldots  \phi (T^{i_n} x)\}$ so 
that $ \phi (T^{i_n} x)=M_n (x)$.  We will this denote this rearrangement by $\{M_0^n (x),  M_1^n (x),\ldots, M_n^n  (x)\}$. Note that $M_n (x)=M^n_n (x)$ in this notation.
Almost sure limit theorems for  trimmed sums involve two sequences of constants
$a(n)$,$b(n)$ so that the scaled truncated sum $\frac{1}{a(n)}\sum_{j=0}^{n-b(n)} M_j^n$ satisfies a strong law of large numbers. We refer especially to~\cite{Aaronson_Nakada,Schindler},
where very precise information on the limiting behavior and choice of constants $a(n)$,$b(n)$ is given for certain dynamical systems. Such results make clear the relations between large extremal values of the time series and the behavior of the Birkhoff sum. There still remains the question of the rate of growth of $\sum_{j=b(n)+1}^{n} M_j^n$.  However
good estimates on the lower bound of the rate of growth of $S_n$ are given by the constants $a(n)$ in the trimmed sum limit. In fact~\cite[Theorem 1.8]{Schindler}  provides a better bound for $\liminf S_n$ in  the
context of piecewise uniformly expanding interval maps than our techniques. We remark on this at more length later. 

In this paper we will consider the observable $\phi (x)=d(x,p)^{-k}$  over chaotic dynamical systems $(T,X,\mu)$  for values of $k$ which ensure that $\int \phi d\mu=\infty$. Our results are limited to probability spaces, in that $\mu (X)=1$. Most of our results
generalize in an obvious way to a wider class of functions, for example those for which $\mu (\phi >t)=\frac{L(t)}{t^{\gamma}}$ where $0<\gamma <1$ and $L(t)$ is a slowly varying function, as long as
the sets $(\phi >t)$ for large $t$ correspond to  sets for which the SBC property holds. Similarly our results generalize to observables $\phi$ 
with a finite set of singularities $\{p_1,\ldots, p_m\}$ such that for all $i$ there
exist constants $C_1$, $C_2$ $r>0$ such that  $0<C_1<\frac{\phi (y)}{d(y,p_i)^{-k}}<C_2$ for all $y\in B(p_i,r)$ and with integrable negative part i.e. if $\phi^{-}:=\max\{ 0,-\phi\}$
then $\int \phi^{-}d\mu<\infty$. But for simplicity of exposition we will stick to the $\phi (x)=d(x,p)^{-k}$.  The case where $\phi^{-}(x)$ is not integrable is very interesting but the techniques of this paper 
are not immediately applicable to this case. We refer to~\cite{Aaronson_Kosloff, Kosloff} for interesting recent results on trimmed symmetric Birkhoff sums in the setting of infinite ergodic theory (when the underlying probability space has infinite measure).

\section{Dynamical Borel Cantelli lemmas and infinite Birkhoff sums.}

We assume that $(T,X,\mu)$ is an ergodic dynamical system and $X$ is a measure and metric space with a Riemannian metric $d$. Let $m$ denote Lebesgue measure on $X$. Let $B(p,r):=\{ x: d(p,x)<r\}$ denote the ball 
of radius $r$ about a point $p$ with respect
to the  given metric $d$. Suppose that  $b(n)$ is a sequence of nested sets in $X$ based about a point $p$.  Define 
\[
E_n=\sum_{j=1}^n \mu (B_j)
\]
For the purposes of this paper (see~\cite{Chernov_Kleinbock}, who introduced the term) we say that the Strong Borel Cantelli (SBC) property holds for  $(B_j)$ if  for $\mu$ a.e. $x\in X$ 
\[
\sum_{j=1}^{n} 1_{B_j}\circ T^j (x)=E_n+ o(E_n)
\]
In most of the examples we consider we have a better  estimate of the error term and, for any $\delta>0$,
\[
\sum_{j=1}^{n} 1_{B_j}\circ T^j (x)=E_n+ O(E_n^{1/2+\delta})~(*)
\]
If $(*)$ holds we say that the sequence $(B_j)$ satisfies the QSBC property, for quantitative Strong Borel Cantelli property. If $T^j (x)\in B_j$ infinitely often for $\mu$ a.e. we say  that the 
sequence $(B_j)$ has the Borel-Cantelli property.

Examples of systems for which the  QSBC property has been proved for balls nested at   points $p$ in phase space  include Axiom A diffeomorphisms~\cite{Chernov_Kleinbock}, uniformly partially hyperbolic systems
preserving a volume measure with exponential decay of correlations~\cite{Dolgopyat}, uniformly expanding $C^2$ maps of the interval~\cite{Phillipp}, and Gibbs-Markov type maps of the interval~\cite{Kim}.
For intermittent type maps with an absolutely continuous invariant probability measure the work of Kim~\cite{Kim} and Gou\"ezel~\cite{Gouezel} gives a fairly complete picture: the Borel-Cantelli
property holds for nested balls except those based at the indifferent fixed point.  Other results on non-uniformly expanding 
systems include one-dimensional maps modeled by Young towers with exponential decay of correlations~\cite{GNO}, the general framework of~\cite{HNPV} and  other hyperbolic settings~\cite{Galatolo,Luzia,Maucourant, Kessebohmer}.

\subsection*{Non-integrable observations.}

Let $\phi (x)=d(x,p)^{-k}$ for some distinguished point $p$, where  $\mbox{dim}(X)=D$ and $k\ge D$.
Let $S_n =\sum_{j=0}^{n-1} \phi \circ T^j$.

\begin{thm}\label{ONE}

Suppose that $(T,X,\mu)$ is an ergodic dynamical system with $\mbox{dim}(X)=D$. Let $\phi (x)=d(x,p)^{-k}$ for some distinguished point $p$.
Suppose  there exist constants $C_1,C_2$ such that $0<C_1<\frac{d\mu}{dm}(p)<C_2$ and that the SBC property holds for nested balls about $p$. 

If $k>D$ then 
 for any $\epsilon>0$
 \[
(a) \limsup \frac{S_n}{n^{k/D}[\log (n)]^{k/D+\epsilon}}=0
\]
 
and for any $\epsilon>0$ 
\[
(b) \liminf \frac{S_n}{n^{k/D-\epsilon}} \ge 1
\]
while 
\[
(c)~S_n \ge n^{k/D} \log^{k/D} n~\mbox{infinitely often}
\]

If moreover the QSBC property holds for nested balls about $p$ then for any $\epsilon>0$ 
\[
(d) \liminf \frac{S_n}{n^{k/D}(e^{-(\log n)^{\frac{1}{2}+\epsilon} })^{k/D}}>1
\]

If $k=D$ the lower bounds in $(b)$ and $(d)$ may be replaced by $\liminf \frac{S_n}{n} \ge >\delta>0$ for some constant $\delta$, while $(a)$ and $(c)$ hold.

\end{thm}

\noindent {\bf Proof.}

We assume first  $k>D$.
It is known from Aaronson~\cite[Proposition 2.3.1]{Aaronson} that if $a(x)$ is increasing, $\lim_{x\to \infty} \frac{a(x)}{x}=0$ and
\[
\int a(\phi (x)) d\mu <\infty
\]
then for $\mu$ a.e. $x$ 
\[
\lim_{n\to \infty} \frac{a(S_n)}{n}=0
\]

Our assumptions imply that $\mu(B(p,r)) \sim r^{D}$. In fact using spherical coordinates our assumption on the density  implies that
$\int d\mu = \int h(x) dx= \int K(\theta_{1},...,\theta_{D-1})r^{D-1} h(r) dr  \theta_1\ldots d\theta_{D-1}$ where $0<c_1 < K(\theta_{1},...,\theta_{D-1})<c_2$ for some constants $c_1$,$c_2$.

 By the Borel Cantelli lemma $\mu (T^n x\in B(p,\frac{1}{n^{1/D+\delta}})~i.~0)=0$ for any $\delta>0$. Hence given $\delta>0$ for $\mu$ a.e. $x\in X$ there exists a time $N(x)$ such that 
$T^{i}x\not \in  B(p,\frac{1}{n^{1/D+\delta}})$ for all $i>N(x)$. This implies that $\phi \circ T^{j} \le n^{k(1/D+\delta)}$ for all $j\ge N(x) $. Thus $S_{n} \le C(x) n^{1+k(1/D+\delta)}$ for large $n$ where $C(x)$ is a constant.  Hence $\log(S_{n}) \le c(x) \log(n)$ for some constant $c(x)>0$. Choosing $a(x)=\frac{x^{D/k}}{\log (x)^{1+\eta}}$ for $\eta>0$  then 

\[
a(S_{n}) = \frac{(S_{n})^{D/k}}{\log(S_{n})^{1+\eta}} \ge \frac{(S_{n})^{D/k}}{[c(x)\log(n)]^{1+\eta}}
\]
Hence 
  for any $\epsilon>0$
 \[
\limsup \frac{S_n}{n^{k/D}[\log (n)]^{k/D+\epsilon}}=0
\]

Assume now that the SBC property holds for nested balls about $p$. 
First note that if $r_n=(n)^{-1/D}$ then $T^n x \in B(p,r_n)$ i.o. Let $B_j:=B(p,\frac{1}{j^{1/D}})$. From the  SBC property 
$\sum_{j=1}^n 1_{B_j} \circ T^j (x) \sim \log (n)$.

If we define  $n_l:=\max_{j\le n} \{ T^j x \in B(r_j,p)\}$ then for $\mu$ a.e. $x\in X$,  for any $M>0$, $\lim{n\to \infty} \frac{n_l}{n^{1-\delta}} >M$ for any $\delta >0$.  To see  this, for a generic $x\in X$, $\lim_{n\to \infty}\frac{S_n}{\log n} =1$. By definition of $n_l (x)$, $S_{n_l}=S_n$ and hence $\lim_{n\to \infty} \frac{S_{n_l}}{\log n}=1$. As $\lim_{n\to \infty} \frac{S_{n_l}}{\log n_l}=1$ we see
$\lim_{n\to \infty} \frac{\log n_l}{ \log n}=1$, which implies the result.

Since $S_n > M_{n_l}$,  $\liminf \frac{S_n}{n^{k/D-\epsilon}} \ge 1$ for any $\epsilon>0$.

 Suppose now that we have a quantitative error estimate in the form of the QSBC property, 
 \[
\sum_{j=1}^{n} 1_{B_j}\circ T^j (x)=E_n+ O(E_n^{1/2+\delta})
\]
Then 
\[
S_n=E_n+ O(E_n^{1/2+\delta})
\]
\[
S_{n_l}=E_{n_l}+ O(E_{n_l}^{1/2+\delta})
\]
By definition of $n_l$, $S_{n_l}=S_n$ and hence
\[
E_n -E_{n_l} =O(E_n^{1/2+\delta})
\]
We obtain
\[
\log n- \log n_l =O(E_n^{1/2+\delta})
\]
which implies that
\[
n_l \ge ne^{-(\log n )^{\frac{1}{2}+\delta}}
\]
for any $\delta>0$.

 Hence $\liminf \frac{S_n}{n^{k/D}(e^{-(\log n)^{\frac{1}{2} +\epsilon} })^{k/D}}>1$ for any $\epsilon>0$.
 
 The proofs of $(a)$ and $(c)$ in the case $k=D$ are unchanged, and estimates  $(b)$ and $(d)$ are immediate  consequences of the ergodic theorem.
 
 \begin{rmk}
 
 The assumptions of Theorem~\ref{ONE} are satisfied by Anosov diffeomorphisms~\cite{Chernov_Kleinbock}, uniformly expanding $C^2$ maps of the interval~\cite{Phillipp} and  Gibbs-Markov type maps of the interval~\cite{Kim}. Kim also shows that for all $p \in X $ in a class of intermittent maps preserving an absolutely continuous probability measure the conditions hold,  except at the indifferent fixed point. Recent work of Tanja Schindler~\cite[Theorem 1.8]{Schindler} on trimmed Birkhoff sums has shown that for Gibbs-Markov maps the limit infimum estimate $d$ can be improved to $\liminf \frac{S_n(\log^k n)}{n^k}>1$.

 \end{rmk}
 
\subsection{Non-integrable observables on a class of intermittent type maps}.

A simple model of intermittency, a form of Manneville-Pommeau map, is the class of maps  $T_{\alpha}$ introduced by  Liverani, Saussol and Vaienti in~\cite{LSV}.
\begin{equation}\label{PM} T_{\alpha}(x)= \begin{cases} x + 2^{\alpha}x^{1+\alpha}, \
0\le x \le 1/2\\ 2x-1, \ 1/2 \le x \le 1
  \end{cases} \qquad 0 \le \alpha <1.
\end{equation}

The map $T_{\alpha}$ has a unique absolutely continuous probability measure $\mu$ if $0\le \alpha <1$. We will only consider the case of a probability measure, rather than an infinite measure preserving system.
The density $h_{\alpha}(x)$  is Lipschitz and strictly positive on any interval of form $[a,1]$, $a>0$ but
blows up at $x=0$, where $h_{\alpha}(x)\sim x^{-\alpha}$.

Kim~\cite[Proposition 4.1]{Kim} has shown that if $p\not =0$ then any nested sequence of balls about $p$ has the SBC property.

\begin{thm}\label{intermittent}
Suppose $(T_{\alpha},[0,1],\mu_{\alpha})$ is a Liverani-Saussol-Vaienti map with $0\le \alpha <1$. Let $p \in [0,1]$ and 
$\phi (x)=d(x,p)^{-k}$ with $k\ge 1$. Define $S_n=\sum_{j=1}^{n}\phi\circ T^j $. Then if $p\not =0$, for any $\epsilon>0$

\[
\liminf \frac{S_n}{n^k(e^{-(\log n)^{\frac{1}{2}+\epsilon}})^{k}} \ge 1
\]
 and
 \[
\limsup \frac{S_n}{n^{k}[\log (n)]^{k+\epsilon}}=0
\]
In particular 
\[
\lim_{n\to \infty} \frac{\log S_n}{\log n}=k
\]
If  $p=0$  then for any $\epsilon>0$
\[
\liminf \frac{S_n}{n^{k+\alpha -\epsilon}}>1
\]
and 
  \[
  \limsup \frac{S_n}{ n^{k+\alpha+\epsilon}}=0
\]
  In particular
\[
\lim_{n\to \infty} \frac{\log S_n}{\log n}=k+\alpha
\]

\end{thm}

\noindent {\bf Proof of theorem.} 

We first consider the case $p\not=0$ and recall a proposition from~\cite{GNO}. We will use it to improve the SBC property  estimate of Kim~\cite[Proposition 4.1]{Kim}
to the QSBC property.

\begin{prop}\label{GM}
  Let $X$ be a compact interval and let $\mathcal{P}$ be a countable partition of $X$ into subintervals.
  Suppose that $(T,X,\mu,\mathcal{P})$ is a Gibbs-Markov system.  Let $(B_{n})$ be a
  sequence of intervals in $X$ for which there exists $C > 0$ such that $\mu (B_{j}) \le C \mu (B_{i})$ for
  all $j \ge i \ge 0$.  If $\sum_{n=0}^{\infty} \mu(B_{n}) = \infty$, then denoting  $E_n=\sum_{j=1}^{n} \mu (B_j)$
  for any $\epsilon >0$, 
  \[
  \sum_{j=1}^n 1_{B_j}\circ T^j (x)= E_n +O(E_n^{1/2+\epsilon})
  \]
  for $\mu $ a.e. $x\in X$.

\end{prop}

A first return time Young Tower $(F,\nu,\Delta)$ may be constructed for this class of intermittent maps with base $\Delta=[1/2,1]$~\cite{LSY}. Every point $p\not =0$ has a unique representation
in  such a first return time Tower, in the sense that there is a unique $t$ such that $F^{-t} (p)\in \Delta$. 
Hence Proposition~\ref{GM} shows that if $p\not =0$ and $(b(n))$ is a sequence of nested sequence of balls based about $p$ then 
\[
  \sum_{j=1}^n 1_{B_j}\circ T^j (x)= E_n +O(E_n^{1/2+\epsilon})
  \]
  for $\mu $ a.e. $x\in X$.

Hence by the proof of Theorem~\ref{ONE}  for $\mu$ a.e. $x$  
\[
\liminf \frac{S_n}{(ne^{-(\log n)^{\frac{1}{2}+\epsilon}})^{k}} \ge 1
\]
 for any $\epsilon>0$, and as a consequence of   Aaronson~\cite[Proposition 2.3.1]{Aaronson}
 for any $\epsilon>0$
 \[
\limsup \frac{S_n}{n^{k}[\log (n)]^{k+\epsilon}}=0
\]

Now we consider the case $p\not =0$.  For nested intervals based at $p=0$ an interesting failure of the  dynamical Borel-Cantelli lemma occurs, described in~\cite{Kim}. 
To understand this phenomenon let $T_1$ and $T_2$ be the two branches of the map
$T_{\alpha}$, with domains $[0,1/2]$ and $[1/2,1]$ respectively. Consider the sequence of sets $b(n)=[0,\frac{1}{n^{\gamma}})$ for any $1<\gamma \le \frac{1}{1-\alpha}$. Kim notes that $\sum_n \mu ( b(n))$
diverges (due to $h_{\alpha}(x)\sim x^{-\alpha}$) while $\sum_n m(b(n))<\infty$.  Note that $T_1^{-1} (b(n)) \subset b(n)$. Hence the only way that $T^j (x)$ can enter $ B_j$ for infinitely many $j$ is that
$T^{j-1} (x)\in T_2^{-1} (B_j)$ for infinitely many $j$. However  the density $h_{\alpha}(x)$ is strictly positive and Lipschitz on any interval $[a,1]$ for $a>0$ and so $\sum_j \mu (T_2^{-1} (B_j))\sim \sum_j m(T_2^{-1}(B_j)) <\infty$ and the 
sequence $(b(n)) $ is not Borel-Cantelli.



 We now consider  the case of $p=0$ and $\phi (x)=d(x,0)^{-k}$. In this setting using Aaronson~\cite[Proposition 2.3.1]{Aaronson} we solve $\int a(\phi)\frac{1}{x^{\alpha}}dx<\infty$ which gives
 an upper bound roughly of form $\limsup \frac{S_n}{n^{k/(1-\alpha)}}=0$, which is not  optimal (being too large as we will see).

To get a better estimate we will  consider the dynamics near the indifferent fixed point.  The following local analysis of a large class of Manneville-Pommeau maps  $T$ (of which the Liverani-Saussol-Vaienti map is a subclass) is taken from~\cite{LSY2}.
 Fix $\epsilon_0>0$, let $x_0\in (0,\epsilon_0]$ and define the sequence $x_n$ by $x_{n-1}=T_{\alpha} x_n$. Young shows that $x_n\sim \frac{1}{n^\beta}$ where $\beta =\frac{1}{\alpha}$. 
 In fact there is a uniform bound on the number of intervals $[\frac{1}{(m+1)^{\beta}},\frac{1}{m^{\beta}}]$ that meet each $[x_{n+1},x_n ]$ and vice-versa.
 
 This implies that if $x=\frac{1}{2}+\frac{1}{2m^{\gamma}}$ then $T_{\alpha}x=\frac{1}{m^{\gamma}}$. Writing $\frac{1}{m^{\gamma}}=x_n$ for some sequence as described above we have $\frac{1}{m^{\gamma}}=\frac{1}{n^\beta}$, hence it takes $n\sim m^{\gamma/\beta}=m^{\gamma\alpha}$ iterates $j$ for $T^{j+1} x$ to escape the region $[0,\epsilon_0]$ i.e. $T^{j+1} x <\epsilon_0 $ for 
 $j <m^{\gamma\alpha}$. Note that $\sum_{j=1}^{n}\phi (x_j)\ge  \sum_{j=1}^{n} j^{k\beta}$ as $x_j\sim \frac{1}{j^{\beta}}$ and hence $S_n\ge  n^{k\beta +1}$. Hence if
  $x=\frac{1}{2}+\frac{1}{2m^{\gamma}}$ then $\sum_{j=1}^{n}\phi\circ T^j x\ge m^{(\gamma/\beta)(k\beta+1)}=m^{\gamma(k+\alpha)}$.
  
  This gives a lower bound on $\liminf S_n$ since if we define $n_l (x) =\max \{ 1\le j\le n: T^j (x)\in [1/2,1/2+\frac{1}{n}]\}$ then  for any $\epsilon>0$, $\liminf\frac{n_l}{n^{1-\epsilon}}\ge1$ by the arguments of the previous section (we use the weaker SBC estimate as the stronger QSBC estimate does not help in this argument).  Furthermore once $T^n x$ enters $[1/2,1/2+\frac{1}{n}]$ it spends $\sim n^{\alpha}$ iterates in the region $(0,\epsilon_0)$ whence $S_{n+n^{\alpha}}\ge n^{k+\alpha-\epsilon}$. As $\alpha <1$ this implies that $\liminf \frac{S_n}{n^{k+\alpha -\epsilon}}>1$ for any $\epsilon>0$.

We will now show $\limsup \frac{S_n}{ n^{k+\alpha+\epsilon}}=0$ for any $\epsilon$, hence $\lim_{n\to \infty} \frac{\log S_n}{\log n}=k+\alpha$.
We first sketch our argument. Let $0<\eta< 1$. Then $\sum_{j=0}^{n-1} 1_{B(p,\frac{1}{j^{\eta}})}\circ T^j(x)\sim n^{1-\eta}$ for $\mu$ a.e. $x$. Note that if $\delta >0$ then  by Borel-Cantelli
$\mu$ a.e. $x\in X$ has the property that $T^n x\in B(p,(n\log^{1+\delta} n)^{-1})$ for only finitely many $n$. Asymptotically almost every $x$ has the property that $T^j x \in B(p,\frac{1}{j^{\eta}})$
 for $\sim n^{1-\eta}$ iterates   $j$  in the interval $1\le j \le n$, after a certain $L(x)$, i.e. for $j\ge L(x)$,  the maximum value that $\phi\circ T^j x$ attains if it enters $B(p,\frac{1}{n^{\alpha}})$ is
$n^k \log^{k(1+\delta)} n$. We break up $S_n$ for large $n$ into the times $j$ that $T^j (x)$ enters $B(p,\frac{1}{n^{\eta}})$, roughly $n^{1-\eta}$ times where the value $\phi\circ T^{j+1} (x)$
is bounded by $n^k \log^{k(1+\delta)} n$ which thus contributes at most $n^{1-\eta}n^{k+\alpha} \log^{k(1+\delta)} n$ to $S_n$  and the times $j$ that $T^j (x)$ enters $B^c(p,\frac{1}{n^{\eta}})$, which
contributes at most $n.n^{\eta(k+\alpha)}=n^{1+\eta(k+\alpha)}$ to the sum $S_n$. Incorporating the log term into the exponent, by choosing $\eta =\frac{k+\alpha}{k+\alpha+1}$ we obtain $\limsup S_n \le n^{k+\frac{1}{k+1}+\alpha}$. 

We will iterate this procedure. Choose $1> \eta_1>\eta_2>\ldots \eta_m>0$ and for simplicity of notation let $B_{\eta_i}=B(p, \frac{1}{n^{\eta_i}})$.

The contribution of the iterates $j$ that enter $B_{\eta_1}$ we bound by the product of the maximum value they may attain, namely  the value $n^{k+\alpha} \log^{k(1+\delta)} n$
and the number of times the point enters this sequence of sets $n^{1-\eta_1}$  to arrive at $n^{k+\alpha+ 1-\eta_1}$ (incorporating  the log term into the exponent). 
This accounts for those iterates that enter $B_{\eta_1} \subset B_{\eta_2}$
and we bound the contribution of those that enter $B_{\eta_2}/B_{\eta_1}$ by $n^{1-\eta_2}.n^{\eta_1(k+\alpha)} =n^{1-\eta_2+\eta_1(k+\alpha)}$. We bound the contribution of those that enter 
$B_{\eta_3}/B_{\eta_2}$ by $n^{1-\eta_3}n^{\eta_2(k+\alpha)}=n^{1-\eta_3+\eta_2(k+\alpha)}$. Continuing in this way we have  a sum of contributions of form $n^{1-\eta_{j+1}+\eta_j(k+\alpha)}$ terminating
with the last contribution, those iterates $j$ that lie in $B_{\eta_m}^c$ whose contribution we bound by $n^{\eta_m(k +\alpha)}.n=n^{1+\eta_m(k+\alpha)}$. 

If $k\ge 1$, choosing $\epsilon=\frac{1}{(k+\alpha)^m}$ and $\eta_i=1-(k+\alpha)^{i-1}\epsilon$ for $i=1,\ldots,m$ the leading term is $n^{k+\alpha+\epsilon}$ corresponding to $n^{k+\alpha+1-\eta_1}$, this $\limsup S_n \le m n^{k+\alpha+\frac{1}{(k+\alpha)^m}}$
which implies the result as $m$ was arbitrary.

\begin{figure}[h!]
	\centering
	\begin{minipage}{.48\textwidth}
		Liverani-Saussol-Vaienti Map.
		\par\medskip
		\includegraphics[width = \textwidth]{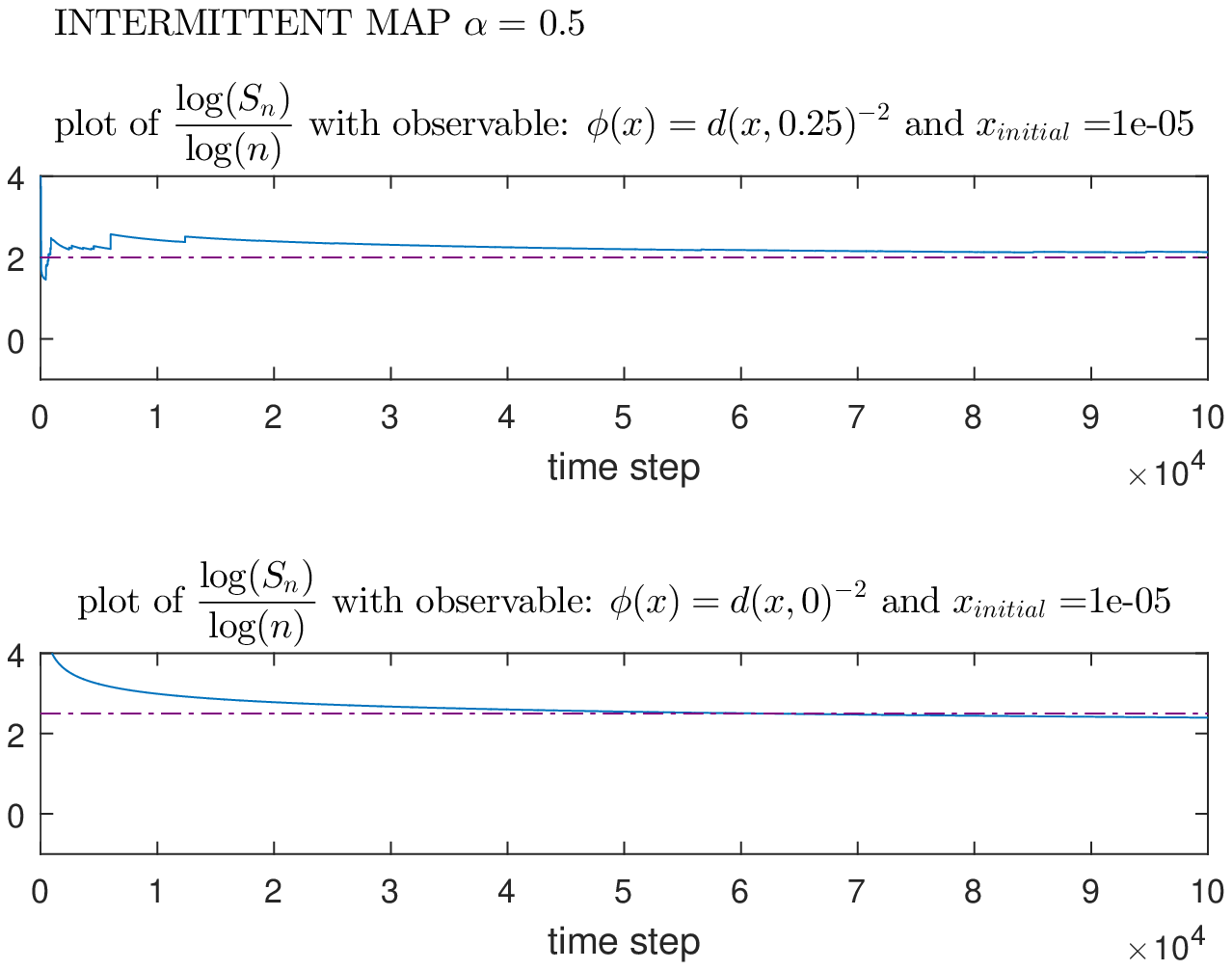}
		\captionsetup{labelformat=empty}
		\caption{Figure: 1(a)}
	\end{minipage}
	\begin{minipage}{.48\textwidth}
		\par\medskip
		\par\medskip
		\par\medskip
		\includegraphics[width = \textwidth]{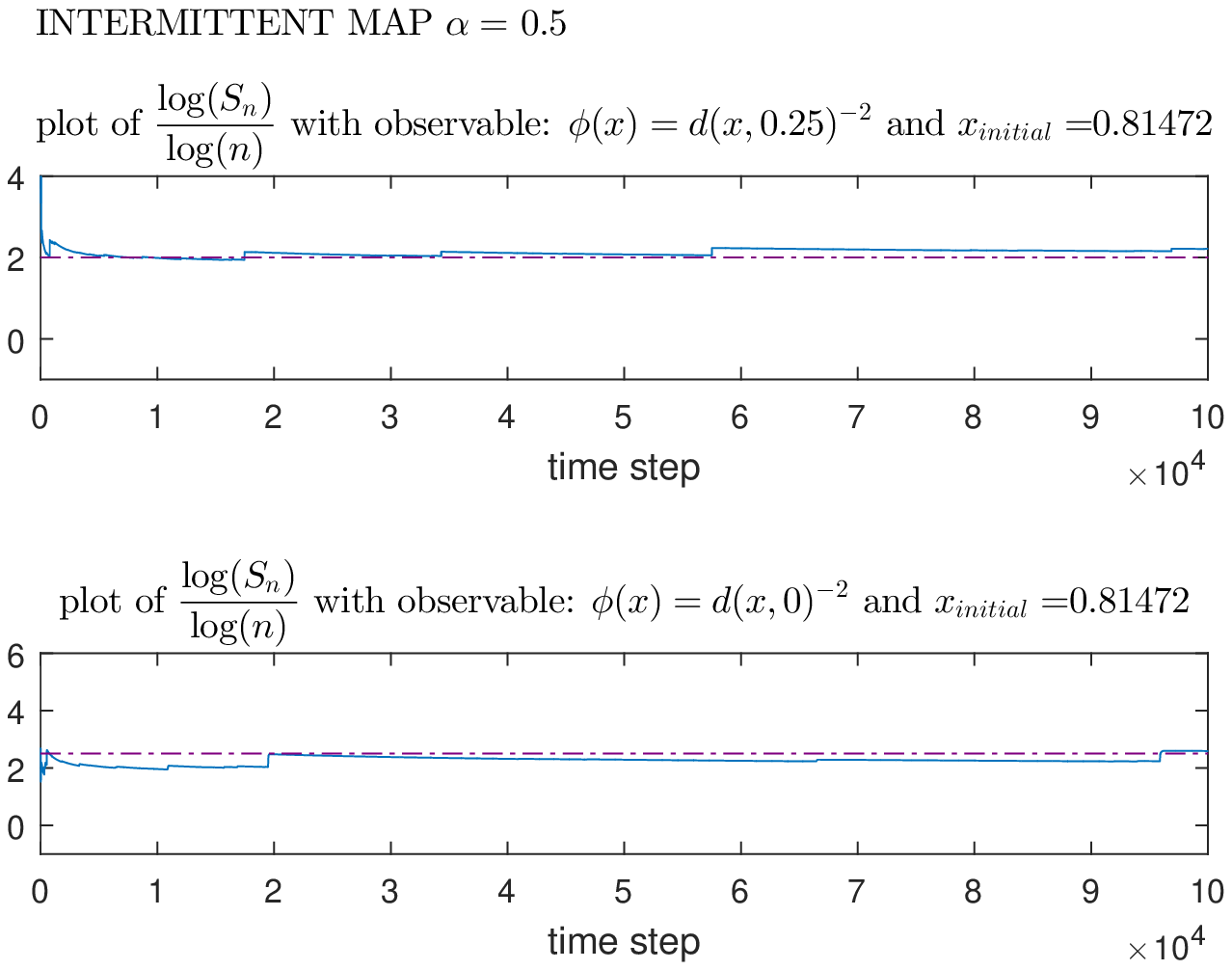}
		\captionsetup{labelformat=empty}
		\caption{Figure: 1(b)}
	\end{minipage}
\end{figure}

\begin{figure}[h!]
	\centering
	\begin{minipage}{.48\textwidth}
		\includegraphics[width = \textwidth]{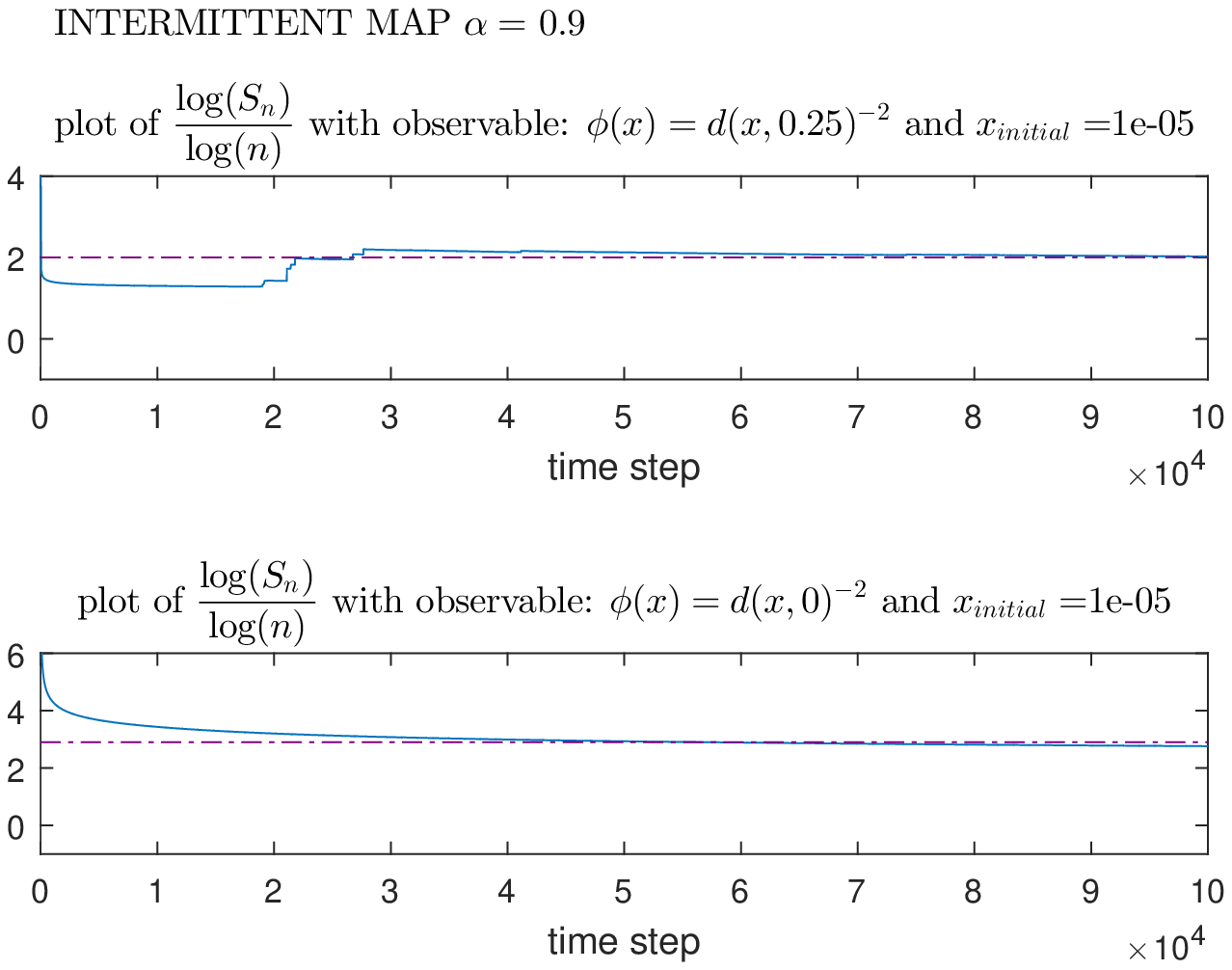}
		\captionsetup{labelformat=empty}
		\caption{Figure: 1(c)}
	\end{minipage}
	\begin{minipage}{.48\textwidth}
		\includegraphics[width = \textwidth]{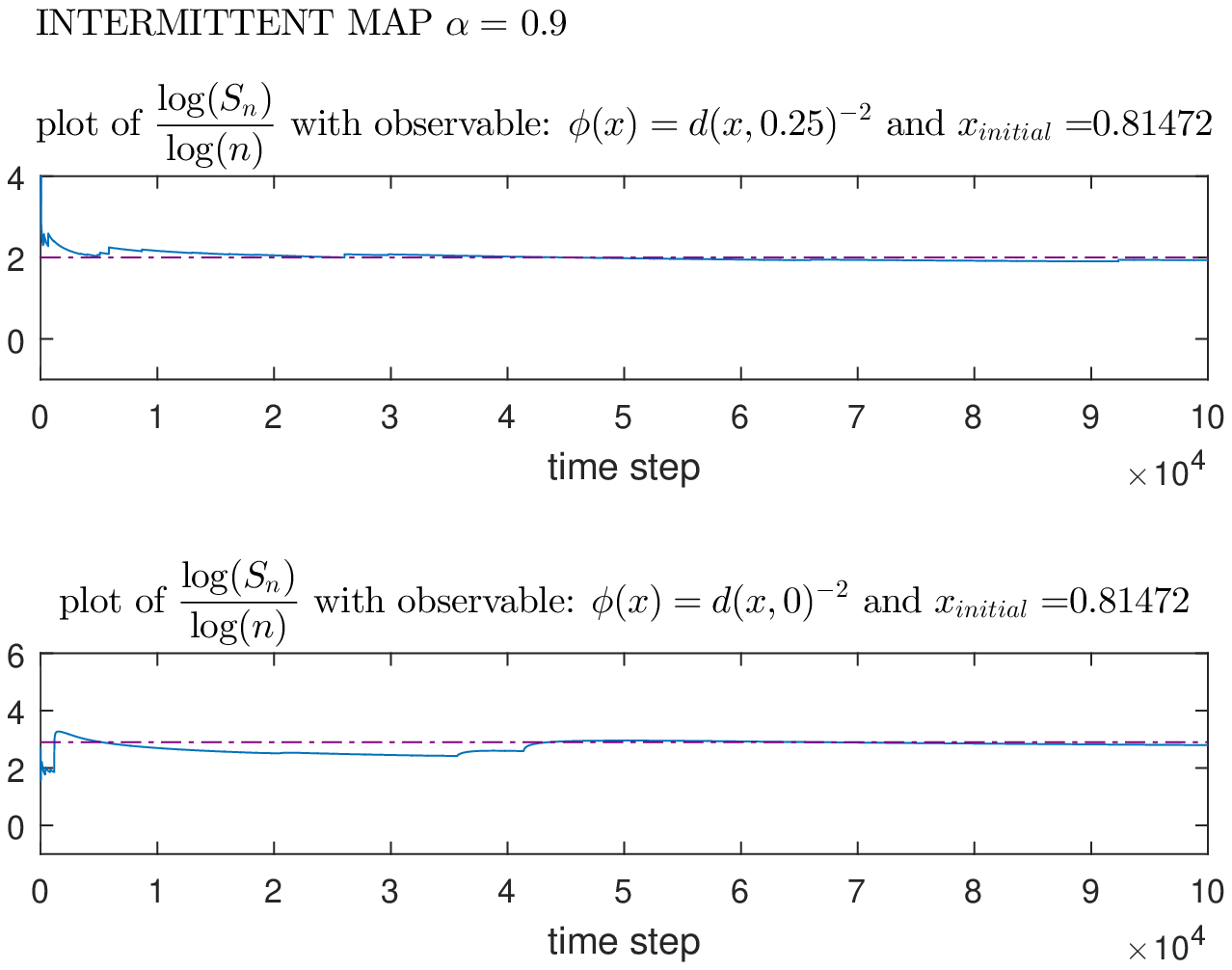}
		\captionsetup{labelformat=empty}
		\caption{Figure: 1(d)}
	\end{minipage}
	\captionsetup{labelformat=empty}
	\caption{Predicted value of the limit is shown by a dotted line.}
\end{figure}

\subsection*{Dynamical systems with $L^p$, $p>1$,  densities and exponential decay of correlations.}

In this section we consider  dynamical systems with exponential decay of correlations, which possess absolutely continuous 
invariant measures (with respect to Lebesgue measure $m$) with densities $\frac{d\mu}{dm}$   in $L^p (m)$.

Suppose $(T,X,\mu)$ is an ergodic measure preserving map of a
probability space $X$ which is also a metric space.  We assume:
\begin{itemize}
\item[(A)] For all Lipschitz functions $\phi,\psi$ on $X$ we have
  exponential  decay of correlations in the sense that there exist
  constants $C,0<\theta<1$ (independent of $\phi$, $\psi$) such that
  \[
  |E( \phi~\psi\circ T^k) -E(\phi) E(\psi) | < C \theta^k\| \phi
  \|_{\mathrm{Lip}} \| \psi \|_{\mathrm{Lip}}.
  \]

\item[(B)] There exist $r_0>0$, $0<\delta <1$ such that for all $p\in
  X$ and all $0<\epsilon< r \le r_0 $
  \[
  \mu\{\, x: r< d(x,p) < r+\epsilon \,\} < \epsilon^{\delta}. 
  \]
\end{itemize}

\begin{lemma}
Suppose $m$ is Lebesgue measure on a $D$-dimensional manifold $X$ and  $h:=\frac{d\mu}{dm} \in L^p (m)$.
Then for all $0<r< r_0$
\[
  \mu\{\, x: r< d(x,p) < r+\epsilon \,\} < \epsilon^{\delta} 
  \]
  for some $\delta>0$

\end{lemma}

\noindent {\bf Proof of lemma:}
Let $q$ be the conjugate of $p$, so that $\frac{1}{p}+\frac{1}{q}=1$. Then $\int_{B_{r+\epsilon}(p)/B_r (p)} d\mu=\int_{B_{r+\epsilon}(p)/B_r (p)} hdx\le 
\|h\|_p m(x: r< d(x,p) < r+\epsilon)^{\frac{1}{q}}$ which implies the result.

Under assumptions $(A)$ and $(B)$  Haydn,  Nicol, Persson  and Vaienti~\cite{HNPV} showed: 

\begin{prop}\label{thm:log}
Assume $(T,X,\mu)$ satisfies assumptions~(A)
and (B).
Suppose $\mu (B_i)\ge C\frac{\log^\beta i}{i}$ for some $\beta>0$, then if  $E_n=\sum_{j=1}^n \mu (B_j)$ for $\mu$ a.e.\ $x \in X$.

\[
\sum_{j=1}^{n} 1_{B_j} \circ T^j (x)=E_n + O(E_n^{1/2+\epsilon})
\]
for any $\epsilon>0$.

\end{prop}
\begin{rmk}
 
 Any exponentially mixing volume preserving system satisfies (A) and (B), for example Sinai dispersing billiard maps with finite and infinite horizon~\cite{LSY,Chernov}.
  Furthermore for a volume preserving dynamical system the density $h(x)=\frac{d\mu}{dm}$
 of the invariant measure is bounded above and is strictly positive. We consider the consequences of this in the next theorem.

  \end{rmk}

 \begin{thm}
 Suppose a dynamical system $(T,X,\mu)$ satisfies (A) and (B) and $p\in X$ has density $h=\frac{d\mu}{dm}$ satisfying $0<C_1<h(p) <C_2$ for some constants $C_1,~C_2$. Suppose also 
  $\mbox{dim}(X)=D$. Then if $\phi (x)=d(x,p)^{-k}$, $k> D$,
  \[
\limsup \frac{S_n}{n^{k/D}[\log (n)]^{\frac{k}{D}+\epsilon}}=0
\]
and 
\[
 \liminf \frac{S_n }{n^{k/D}(e^{-[\log (n)]^{\frac{1}{2}+\epsilon} }  )^{\frac{k}{D}} }>1
 \]
  for any $\epsilon>0$.
  
  In the case $k=D$
  \[
 \liminf \frac{S_n }{n }>a>0
 \]
 for some $a>0$.

 \end{thm} 
 
 \begin{rmk}
 
 By ergodicity  in the case $k=D$
  \[
 \liminf \frac{S_n }{n }>a>0
 \]
 for some $a>0$.

 \end{rmk}

 \noindent {\bf Proof.}
 
 Note that  if $B_r$ is a ball of small radius $r>0$ nested at $p$ then $\mu (B_r)\sim C r^{D}$. We first consider the case $k>D$.
 Let $\phi (x)=d(x,p)^{-k}$ and  $a(x)=|x|^{\frac{D}{k}}/(|\log|x||)^{1+\eta}$. Then $\int a(\phi (x))dx<\infty$. If we define 
$S_n=\sum_{j=1}^n \phi \circ T^j$, then by~\cite[Proposition 2.3.1]{Aaronson} 
\[
\frac{a(S_n)}{n}\to 0
\]
for $\mu$ a.e. $x\in X$. Hence for any $\epsilon>0$, for $\mu$ a.e. $x \in X$
 \[
\limsup \frac{S_n}{n^{k/D}[\log (n)]^{\frac{k}{D}+\epsilon}}=0
\]


To obtain a limit infimum estimate we modify our previous argument. Let $b(n)$ be balls of $\mu$ (hence $m$)  measure $\sim \frac{\log^{\beta} n}{n}$ nested about $p$. Let $E_n:=\sum_{j=1}^n \mu (B_j)$

  Define $n_l:=\max \{0\le j\le n: T^j (x) \in B_j \}$ as before we have 

\[
\sum_{j=1}^n 1_{B_j}\circ T^j=E_n+ O(E_n^{1/2+\delta})
\]
\[
\sum_{j=1}^{n_l} 1_{B_j}\circ T^j=E_{n_l}+ O(E_{n_l}^{1/2+\delta})
\]
By definition of $n_l$, $\sum_{j=1}^{n_l} 1_{B_j}\circ T^j=\sum_{j=1}^n 1_{B_j}\circ T^j$ and hence
\[
E_n -E_{n_l} =O(E_n^{1/2+\delta})
\]
We obtain
\[
\log^{1+\beta} n- \log^{1+\beta} n_l =O(\log ^{1/2+\gamma}(n))
\]
where $\gamma=\delta+\frac{\beta}{2}$.
As $x-y \le x^{1+\beta} -y^{1+\beta}$ for large $y$ and $x>y$ 
we see that 
\[
n_l \ge n e^{-(\log n)^{\frac{1}{2}+\epsilon}}
\]
for any $\epsilon>0$.

Note  that balls of radius $r$ based at $p$ satisfy $\mu(B_r(p))\sim C r^{D}$, and so we are able to bound $S_n$ below by $ M_{n_l}\ge (ne^{-[\log (n)]^{\frac{1}{2}+\epsilon} }  )^{\frac{k}{D}} $.

 Hence \[
 \liminf \frac{S_n }{(ne^{-[\log (n)]^{\frac{1}{2}+\epsilon} }  )^{\frac{k}{D}} }>1
 \]
  for any $\epsilon>0$.

 A recent result of J.Rivera-Letelier~\cite[Corollary B]{Letelier_2012} states:
 
 \begin{prop}
 Let $T$  be a non-degenerate smooth interval map having an exponentially mixing absolutely continuous invariant probability measure
  $\mu$. Then there is $p>1$ such that the density $h$  of $\mu$  with respect to Lebesgue measure $m$ is in $L^p (m)$.
  Moreover, $\mu$  can be obtained through a Young tower with an exponential tail estimate. In particular, $\mu$  satisfies the local central limit theorem and the vector-valued almost sure invariance principle.
  \end{prop}
 
For such maps if the invariant density at   $p$ satisfies $h(x)\sim C d(p,x)^{-\alpha}$, $\alpha >0$, then we have  the estimates:

\begin{thm}\label{density}
 Suppose a dynamical system $(T,X,\mu)$ satisfies (A) and (B) and $p\in X$ has density satisfying $h(x)\sim C d(p,x)^{-\alpha}$, $\alpha >0$. Suppose also 
  $\mbox{dim}(X)=D$. Then if $\phi (x)=d(x,p)^{-k}$, $k\ge D-\alpha$, 
  \[
\limsup \frac{S_n}{n^{k/(D-\alpha)}[\log (n)]^{k+\epsilon}}=0
\]
and 
\[
 \liminf \frac{S_n }{n^{k/(D-\alpha)}(e^{-[\log (n)]^{\frac{1}{2}+\epsilon} }  )^{\frac{k}{D-\alpha}} }>1
 \]
  for any $\epsilon>0$.
  Hence 
  \[
  \lim{n\to \infty} \frac{\log S_n}{\log n}=\frac{k}{D-\alpha}
  \]

 \end{thm} 
  \begin{rmk}
 
 By ergodicity  in the case $k=D$
  \[
 \liminf \frac{S_n }{n }>a>0
 \]
 for some $a>0$.

 \end{rmk}
 
 \begin{rmk}
 This result contrasts with that of the intermittent map where at the indifferent fixed point $x=0$, with density $h(x)\sim x^{-\alpha}$ it was shown  
 for the observable $\phi (x)=x^{-k}$ that
 $\lim_{n\to \infty} \frac{\log S_n}{\log n}=k+\alpha$.
 \end{rmk}

 \noindent {\bf Proof:}
 
 The proof is an obvious modification of the proof of the previous theorem. Let $\tilde{D}=D-\alpha$ and 
 define $a(x)=\frac{|x|^{\tilde{D}/k}}{(\log |x|)^{1+\eta}}$. Then $\int a(\phi (x)) dx<\infty$ and by ~\cite[Proposition 2.3.1]{Aaronson} 
 $\frac{a(S_n)}{n}\to 0$ and hence 
  \[
\limsup \frac{S_n}{n^{k/(\tilde{D})}[\log (n)]^{k+\epsilon}}=0
\]

We now obtain our limit infimum estimate. 

Let $b(n)$ be balls of $\mu$ measure 
$\sim \frac{\log^{\beta} n}{n}$ nested about $p$.  Define $n_l:=\max \{0\le j\le n: T^j (x) \in B_j \}$ as before we have 

\[
\sum_{j=1}^n 1_{B_j}\circ T^j=E_n+ O(E_n^{1/2+\delta})
\]
\[
\sum_{j=1}^{n_l} 1_{B_j}\circ T^j=E_{n_l}+ O(E_{n_l}^{1/2+\delta})
\]
and hence
\[
E_n -E_{n_l} =O(E_n^{1/2+\delta})
\]
We have
\[
\log^{1+\beta} n- \log^{1+\beta} n_l =O(\log ^{1/2+\gamma}(n))
\]
where $\gamma=\delta+\frac{\beta}{2}$.
As $x-y \le x^{1+\beta} -y^{1+\beta}$ for large $y$ large and $x>y$ 
we see that as in the previous theorem
\[
n_l \ge n e^{-(\log n)^{\frac{1}{2}+\epsilon}}
\]
for any $\epsilon>0$.

Note  that balls of radius $r$ based at $p$ satisfy $\mu(B_r(p))\sim C r^{\tilde{D}}$ we see that $S_n\ge M_{n_l}$ implies
\[
 \liminf \frac{S_n }{(ne^{-[\log (n)]^{\frac{1}{2}+\epsilon} }  )^{\frac{k}{D-\alpha}} }>1
 \]
  for any $\epsilon>0$.

 \begin{figure}[h!]
 	\centering
 	\begin{minipage}{.48\textwidth}
 		\includegraphics[width = \textwidth]{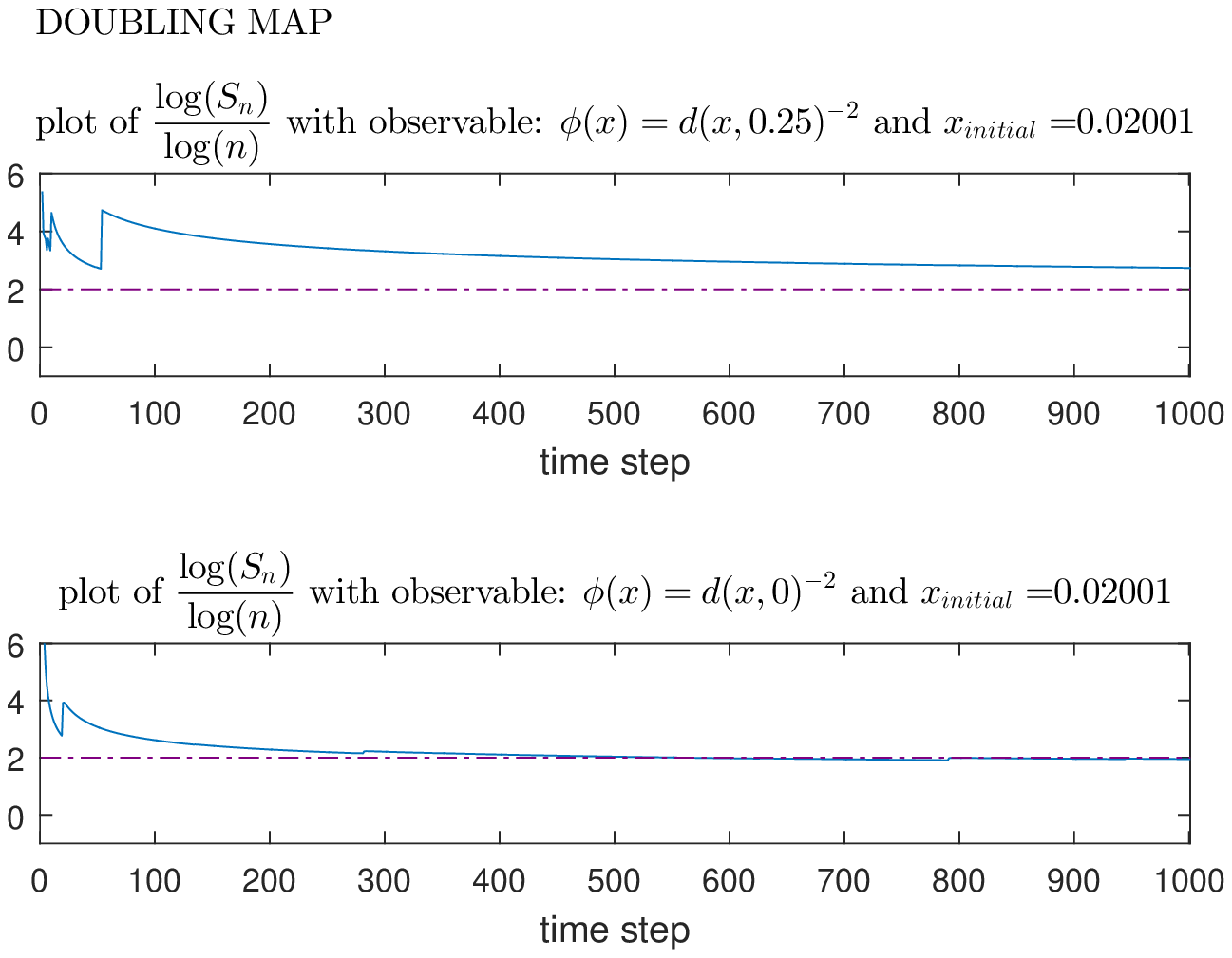}
 		\captionsetup{labelformat=empty}
 		\caption{Figure: 2 Doubling Map.}
 	\end{minipage}
 	\begin{minipage}{.48\textwidth}
 		\includegraphics[width = \textwidth]{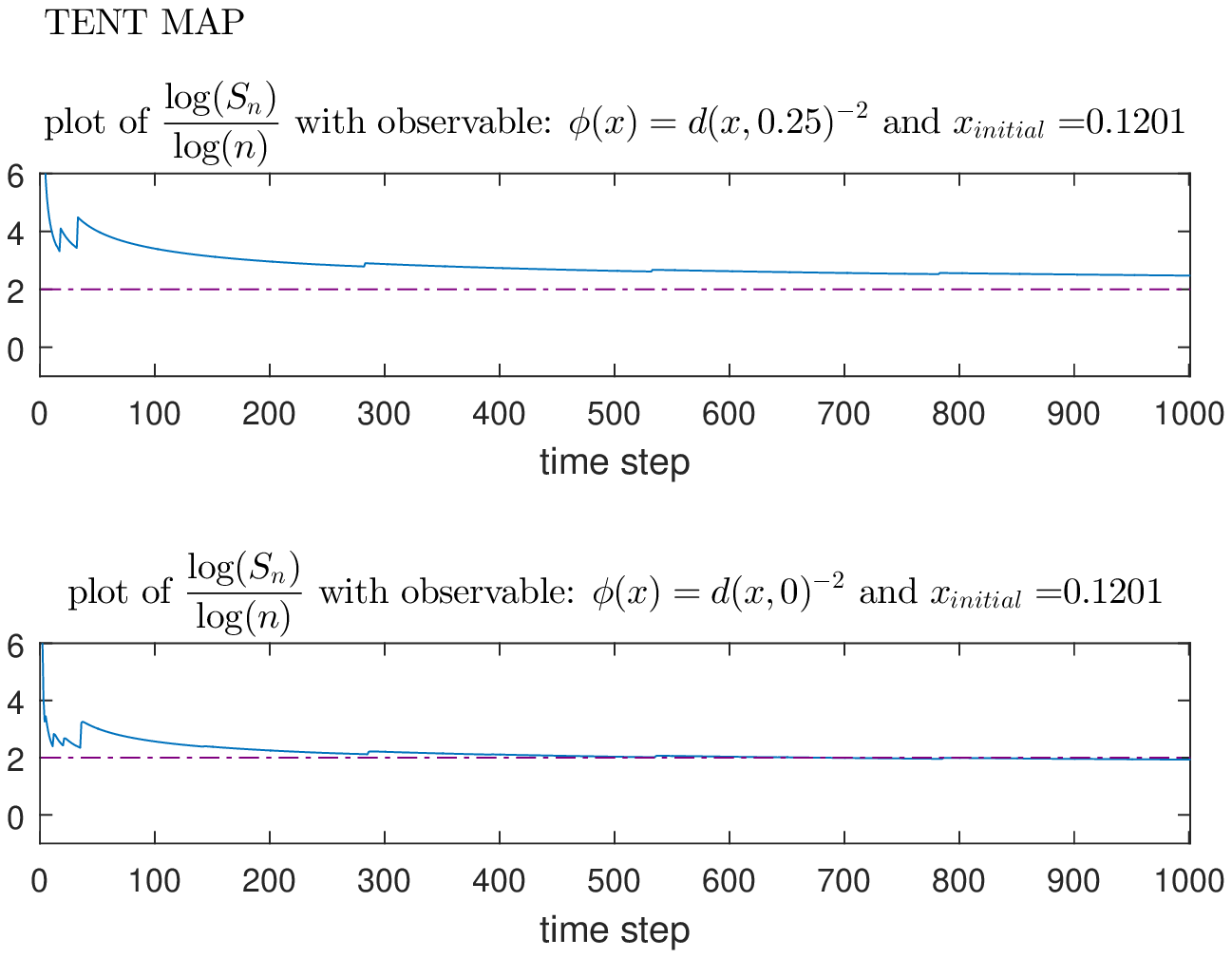}
 		\captionsetup{labelformat=empty}
 		\caption{Figure: 3 Tent Map.}
 	\end{minipage}
 \captionsetup{labelformat=empty}
 \caption{Predicted value of the convergence is marked in a dotted line.}
 \end{figure}
  
\begin{cor}
 Suppose $T(x)=4x(1-x)$ is a unimodal map of the interval $[0,1]$. Let $\phi (x)=d(x,p)^{-k}$,  then if $p=0$ or $p=1$
 \[
 \lim_{n\to \infty} \frac{\log(S_n)}{\log n}=4
 \]
 while  if $p\in (0,1)$

 \[
 \lim_{n\to \infty} \frac{\log(S_n)}{\log n}=2
 \]
 
 \end{cor}
 
   \begin{figure}[h!]
   	\centering
   	\begin{minipage}{.48\textwidth}
   		\includegraphics[width = \textwidth]{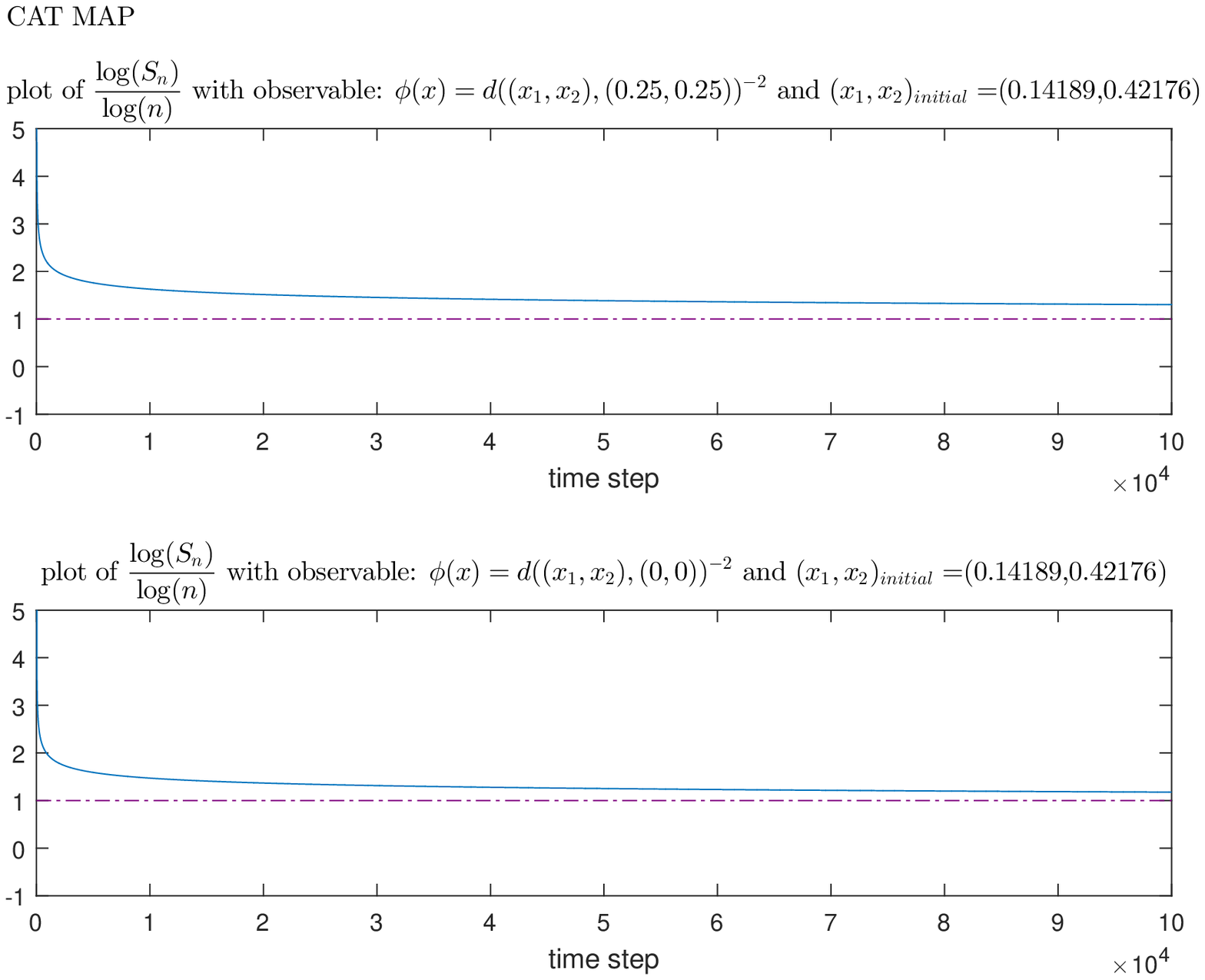}
   		\captionsetup{labelformat=empty}
   		\caption{Figure: 4 Hyperbolic Toral Automorphism.}
   	\end{minipage}
   	\begin{minipage}{.48\textwidth}
   		\includegraphics[width = \textwidth]{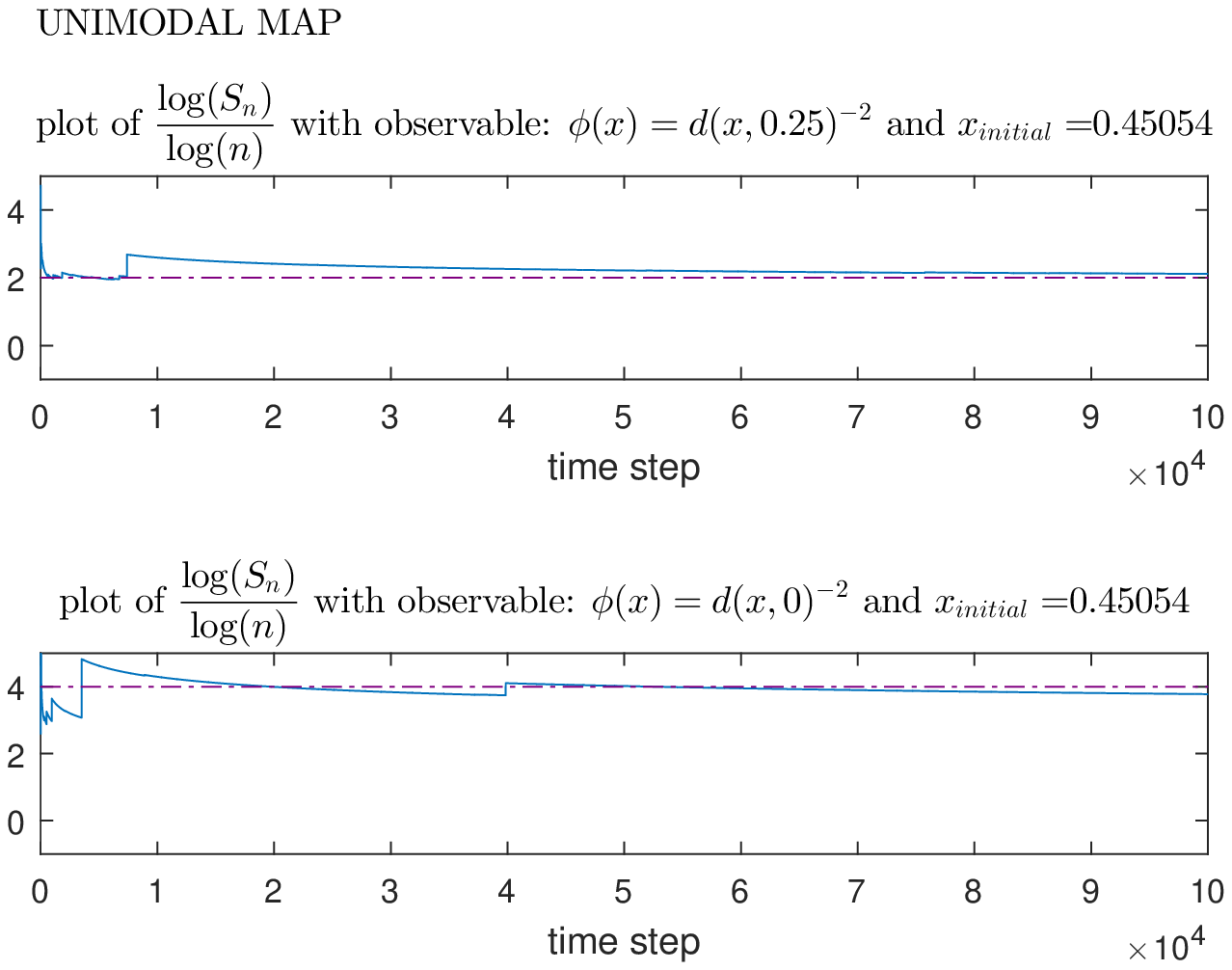}
   		\captionsetup{labelformat=empty}
   		\caption{Figure: 5 Unimodal Map.}
   	\end{minipage}
   	\captionsetup{labelformat=empty}
   	\caption{Predicted value of the convergence is marked in a dotted line.}
   \end{figure}
 
 \noindent {\bf Proof of corollary:}
 
 This map has invariant density $h(x)=\frac{1}{\sqrt{\pi x(1-x)}}$. First note that the unimodal map has density $h(x)\sim\frac{1}{\sqrt{x}}$ for $p = 0$ and $p = 1$. 
which implies the result.\\

 \newpage
 
 



\section{Conclusion.}
 Dynamical Borel Cantelli lemmas and Aaronson~\cite[Proposition 2.3.1]{Aaronson} give useful  bounds on the rate of growth of positive non-integrable functions on
 ergodic dynamical systems. In the case of Gibbs-Markov maps the lower bounds we obtain are  not optimal~\cite{Schindler}. Quantitive Borel-Cantelli estimates and the density of the invariant measure both play a role, for example the contrasting behavior in Theorem~\ref{intermittent} and Theorem~\ref{density}. It would be of interest to develop insights into a broader class
 of examples. It would also be interesting to 
 explore more examples in the setting of Birkhoff sums of  functions $\phi=\phi^{+}-\phi_{-}$ on ergodic probability measure preserving systems with non-integrable positive and
non-integrable  negative parts.

\end{document}